\documentclass[review,onefignum,onetabnum, final]{article}

\usepackage{amsmath,amssymb,graphicx,subcaption,amsbsy}
\usepackage[disable]{todonotes}

\usepackage[numbers]{natbib}
\let\oldref\ref
\renewcommand{\ref}[1]{(\oldref{#1})}

\newcommand{\abs}[1]{\left \lvert #1 \right \rvert}
\newcommand{\dphi}[2]{\phi^{#1}_{#2}}
\newcommand{\ext}[1]{\underset{#1}{\text{ext}}}
\graphicspath{{/figures/}}

\newcommand{\norm}[1]{\left| \left| #1 \right| \right|}

\newcommand{\R}{\mathbb{R}}

\newcommand*{\defeq}{\mathrel{\vcenter{\baselineskip0.5ex \lineskiplimit0pt
                     \hbox{\scriptsize.}\hbox{\scriptsize.}}}%
     		     =}

\def\XXint#1#2#3{{\setbox0=\hbox{$#1{#2#3}{\int}$ }
\vcenter{\hbox{$#2#3$ }}\kern-.6\wd0}}

\title{Optimal Human Navigation in Steep Terrain:\\a Hamilton-Jacobi-Bellman Approach\thanks{Submitted: \today}}

\author{Christian Parkinson\thanks{Department of Mathematics, UCLA, Los Angeles, CA 90095, (chparkin@math.ucla.edu).}%\email{...}
\and David Arnold\thanks{Department of Mathematics, UCLA, Los Angeles, CA 90095, (darnold@math.ucla.edu).}
\and Andrea L. Bertozzi\thanks{Department of Mathematics, UCLA, Los Angeles, CA 90095, (bertozzi@math.ucla.edu).}
\and Yat Tin Chow\thanks{Department of Mathematics, UCLA, Los Angeles, CA 90095, (ytchow@math.ucla.edu).}
\and Stanley Osher\thanks{Department of Mathematics, UCLA, Los Angeles, CA 90095, (sjo@math.ucla.edu).}}

\begin{document}

\maketitle

\begin{abstract}
We present a method for determining optimal walking paths in steep terrain using the level set method and an optimal control formulation. By viewing the walking direction as a control variable, we can determine the optimal control by solving a Hamilton-Jacobi-Bellman equation. We then calculate the optimal walking path by solving an ordinary differential equation. We demonstrate the effectiveness of our method by  computing optimal paths which travel throughout mountainous regions of Yosemite National Park. We include details regarding the numerical implementation of our model and address a specific application of a law enforcement agency patrolling a nationally protected area. 
\end{abstract}

% REQUIRED

%\listoftodos

\section{Introduction} \label{intro}

We consider the problem of determining the optimal walking path between two points given elevation data in a region. If the terrain is fairly flat, this may be very easy as conventional wisdom (``the shortest path between two points is a straight line") will provide a good approximation to the optimal path. However,  in mountainous regions, straight line travel is often inefficient or impossible and the optimal path between two points is no longer clear. 

The problem of optimal path planning goes back at least as far as Dijkstra \cite{Dijkstra} who designed an algorithm for optimally traversing weighted graphs. In the years since, significant effort has been devoted to developing and improving algorithms which find optimal or near-optimal paths in a discrete setting \cite{Hirsch, JSBMitchell, Papadimitrou}. Others have used modified versions of Dijkstra's algorithm for path planning in a semi-continuous setting \cite{IanMitchell, Tsitsiklis}.

Recently, path planning problems have been largely reframed using control theory and partial differential equations.  An early approach was to compute geodesics on triangulated manifolds using an Eikonal equation and gradient descent \cite{KimmelSethian2}. One interesting application of path planning problems is to simple, autonomous robots. These so-called \emph{Dubins' cars} were constrained by a maximum turning radius so Dubins considered paths with bounded local curvature \cite{Dubins}. Recently, this problem was reformulated using a Hamilton-Jacobi-Bellman equation and adapted to include impassable obstacles \cite{AgarwalWang, Takei}. Indeed, Hamilton-Jacobi-Bellman (HJB) equations are now used extensively in optimal path problems. Recent research employs HJB equations in determining reachable and avoidable sets when traveling from a given ground state \cite{TomlinAircraft, Lygeros}. Tomlin et al. also use HJB equations in adversarial reach-avoid games wherein a group of attackers attempt to reach a target set while also avoiding defenders \cite{TomlinReachAvoid}. Others have considered optimal travel in regions which randomly switch between different states; for example, this randomness could account for the affect or weather patterns on a sailboat \cite{ShenVladimirsky}. 

There has been some research into path planning in a geographical or terrain-based setting though most previous work is focused on discrete, graph-based methods employing Djikstra's algorithm and its many variants: so-called $A^*$ and $D^*$ algorithms \cite{krishnaswamy1995resolution, SARANYA2016178}. Such methods have long been used for vehicular navigation and can be adapted to include real-time obstacle recognition \cite{RealTime}. This problem is also of particular interest to those working on UAVs and other autonomous robots \cite{MultiGeo, hachour2008path, planetary, lin2009uav}. In a continuous approach, Popovi\'c et al. \cite{Popovi2017MultiresolutionMA} propose a path-planning algorithm for UAVs by maximizing an information functional which measures the amount of data a UAV can collect. However, the methods of control theory and HJB equations have yet to be applied to terrain-based path-planning which means, for example, that no previous approach has been able to dynamically account for the optimal direction of travel along a path.

In section \ref{model}, we present a model which uses the level set method and a HJB formulation to compute optimal walking paths in a continuous setting where travel direction can be considered dynamically and walking speed is dependent on slope of the local terrain and. This is as opposed to other terrain-based path planning methods which are fully or partially discrete and do not account for directional movement. In section \ref{numerics}, we discuss the numerical simulation of the model. We begin by testing the model against toy problems using synthetic elevation data specifically designed so that the ``correct answer" is somewhat clear \emph{a priori} and move on to use real elevation data of Yosemite National park. Results of numerical simulations are presented in section \ref{results}. The motivation for this work was to aid law enforcement agencies in efficiently patroling protected areas such as parks or forests, but with small adjustments, our method could be applied to optimal path planning in any number of scenarios.

\section{Mathematical Model}\label{model}

Our primary mathematical tool is the level set method of Osher and Sethian \cite{OsherSethian}.  The level set method models propagation of fronts by treating them as the zero level set of auxiliary function $\phi$, known as the \emph{level set function}. We will discuss the method in two spatial dimensions since this is relevant for terrain-based path-planning, but this can be easily generalized to higher dimensions. 

\subsection{The Level Set Method}

Suppose that $\Omega \subset \R^2$ is open and bounded with Lipschitz continuous boundary $\Gamma(0) = \partial \Omega$ which is the curve that will evolve via some level set motion. To begin, we find a Lipschitz continuous function $\phi_0 : \R^2 \to \R$ such that $\phi_0 < 0$ in $\Omega$ and $\phi_0 > 0$ in $\R^2 \setminus \Omega$. Continuity of $\phi_0$ implies that $\Gamma(0) = \{x \in \R^2 : \phi_0(x) = 0\}$; that is, $\Gamma$ is the zero level contour of the initial function. Next, we evolve the function $\phi : \R^2 \times \{t > 0\} \to \R$ using the Hamilton-Jacobi equation \begin{equation} \label{eq:levelSet} \begin{aligned}
\phi_t + H(x,\nabla \phi) = 0,&\\
\phi(x,0) = \phi_0(x),&
\end{aligned}
\end{equation} where the \emph{Hamiltonian} $H(x,p)$ is homogeneous of degree one in the variable $p$; here $p$ is a proxy for $\nabla \phi$. As $\phi$ evolves according to the PDE, we define $\Gamma(t) = \{x \in \R^2: \phi(x,t) = 0\}$ (so that, in particular, $\Gamma(0) = \Gamma$) and the curve $\Gamma(t)$ evolves with level set motion which is prescribed by the Hamiltonian $H$ \cite{osher2003level}.

In the simplest case $H(x,p) = \abs{p}$ and \ref{eq:levelSet} is the Eikonal equation $\phi_t + \abs{\nabla \phi} = 0$. Re-writing the equation as $\phi_t + \hat n \cdot \nabla \phi = 0$ where $\hat n = \nabla \phi / \abs{\nabla \phi}$, it is clear that locally this equation gives advection in the outward normal direction with velocity 1. This causes $\Gamma(t)$ to deform outward with normal velocity 1. In this case, for $t > 0$, $\Gamma(t)$ represents the set of all points which which are distance $t$ from the original curve $\Gamma(0)$. Equivalently, since we are considering people traveling throughout regions, $\Gamma(t)$ is the set of points which can be reached if one travels from $\Gamma(0)$ with normal velocity $1$ for time $t$. To prescribe a different normal velocity $v(x)$ rather than allowing individuals to travel with normal velocity 1, one can simply modify the Hamiltonian by setting $H(x,p) = v(x) \abs{p}$. Now $\Gamma(t)$ represents the set of points which can be reached if one travels from $\Gamma(0)$ with normal velocity $v(x)$ for time $t$.

Using the level set equation, one can compute the (approximate) time that it takes to travel from one point to another in our domain. Let $a \in \R^2$ represent a starting point and $b \in \R^2$ represent an ending point. For some small $\delta > 0$, let $\phi_0(x) = \abs{x-a} - \delta$ so that $\Gamma(0) = \{\abs{x-a} = \delta\}$ is a small circle around the point $a$. When $\Gamma(t)$ evolves outward with prescribed normal velocity $v(x)$, there will be some time $t^* > 0$ such that $b \in \Gamma(t^*)$; that is, at some positive time $t^*$, the level set will hit our ending point. This time $t^*$ is the time required to travel from point $a$ to point $b$ when traveling in the normal direction with velocity $v(x)$ (neglecting the small parameter $\delta$). This gives a method for calculating travel times, but this model is too simple for our purposes, only allowing for travel in the normal direction which is potentially far from optimal. For example, if in a physical setting there is a large mountain between the points $a$ and $b$, one may wish to walk around the mountain rather than over the mountain, as normal direction travel may suggest.  Thus at each point, one must not only consider the speed of travel, but also the direction of travel. Considering direction, it no longer makes sense to simply specify a velocity $v(x)$ at each point. Instead, we assume that walking velocity depends on both the gradient of the terrain at the current point and the direction of travel as we search for the optimal travel direction.

\begin{figure}[b!]
\centering
\includegraphics[width=0.63\textwidth]{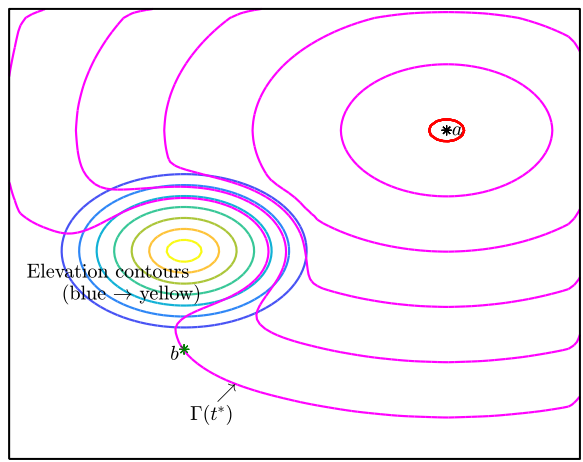}
\caption{To find optimal travel time, begin with a small circle around $a$ (red) and evolve level sets $\Gamma(t)$ (magenta) outward until the time $t^* > 0$ at which $b \in \Gamma(t^*)$.}
\label{key1}
\end{figure}

\subsection{Our Model} \label{ourmodel} For our model, assume that in addition to the starting and ending points $a,b \in \R^2$, there is an elevation profile $E(x)$ and a velocity function $V(S)$ which gives human walking velocity as a function of terrain slope $S$. Let $\theta \in [0,2\pi]$ be a control variable which represents walking direction and let $s(\theta) = (\cos(\theta), \sin(\theta))$ be the corresponding direction vector. Now if one is standing at a point $x$ and desires to walk in the direction $\theta$, they can walk with velocity $V(s(\theta)\cdot\nabla E(x))$ since $s(\theta) \cdot \nabla E(x)$ represents the slope at $x$ in the direction of $\theta$. For each $\theta \in [0,2\pi]$, define the directional Hamiltonian $H_\theta(x,p) = V(s(\theta)\cdot\nabla E(x))[s(\theta) \cdot p]$. Note that using this Hamiltonian, the corresponding Hamilton-Jacobi equation models advection in the direction of $\theta$. To consider optimal travel, we take the supremum over all $\theta$. Define the \emph{optimal path Hamiltonian} \begin{equation} \label{eq:hamiltonian}
H(x,p) \defeq \sup_{\theta \in [0,2\pi]} H_\theta(x,p) = \sup_{\theta \in [0,2\pi]} \big\{ V(s(\theta)\cdot\nabla E(x))[s(\theta) \cdot p] \big\}.
\end{equation} This results in a Hamilton-Jacobi-Bellman equation, which, after resolving the supremum in $\theta$ takes the form \ref{eq:levelSet}; it is indeed a level set equation, since this optimal path Hamiltonian is still homogeneous of degree one in the variable $p$. Now to find the optimal travel time between points $a$ and $b$, one can use the same method described above: letting $\Gamma(0) = \{x \in \R^2 : \abs{x-a} = \delta\}$ for small $\delta$, evolve $\Gamma(t)$ using the level set equation with the optimal path Hamiltonian until the time $t^* >0$ such that $b \in \Gamma(t^*)$. This $t^*$ is the minimal time required to travel from $a$ to $b$. This procedure is displayed in Figure~\ref{key1}. 

What remains is to compute the optimal path from $a$ to $b$: the path which requires time $t^*$ to traverse. In order to do this, one simply needs to follow the characteristics of the Hamilton-Jacobi-Bellman equation. We would like to travel along characteristics originating from $a$ toward $b$. However, with this small parameter $\delta$, we have removed a small neighborhood of $a$ and instead initiate the motion from the circle $\Gamma(0)$. Note for example, that $\phi_0$ is non-differentiable at $a$. Accordingly, one should follow the characteristics backwards from $b$ to $\Gamma(0)$. The characteristic equations are \begin{equation} \label{eq:ODE} \begin{split}\dot x &= -\nabla_p H(x,p), \,\,\,\, x(0) = b, \\ \dot p &=\hphantom{-}\nabla_x H(x,p), \,\,\,\, p(0) = \nabla \phi(b,t^*). \end{split} \end{equation} Physically, one can imagine starting at the point $b$, considering what was the direction of the optimal step which led to the current point, stepping backwards in that direction and updating the direction in real time as one is walking backwards. Running this system of ODEs to time $t^*$, one will have backtracked optimally all the way from $b$ to $\Gamma(0)$.  

To summarize, once we have defined the optimal path Hamiltonian \ref{eq:hamiltonian}, the algorithm for finding the optimal path consists of two steps: \begin{itemize}
\item[1.] Find the optimal travel time by advancing the PDE \begin{align*} \phi_t + H(x,\nabla \phi) = 0,& \\\phi(x,0) = \abs{x-a} - \delta,& \end{align*} until the time $t^* > 0$ such that $b \in \Gamma(t^*)$. \\

\item[2.] Find the optimal travel path by advancing the ODE system \begin{align*} \dot x &= -\nabla_p H(x,p), \,\,\,\, x(0) = b, \\ \dot p &=\hphantom{-}\nabla_x H(x,p), \,\,\,\, p(0) = \nabla \phi(b,t^*) \end{align*} until time $t^*$.
\end{itemize} 

\subsection{The Associated Control Problem} Since we are determining optimal travel, we know that underlying the formalism of section \ref{ourmodel}, there is a control problem that is being solved and a payoff function which is being maximized. As above, let $a \in \R^2$ be the initial point. If one is standing at the point $x \in \R^2$, then traveling optimally away from the point $a$ for a time $t$ is the same as maximizing the distance $\abs{\mathbf{x}(t) - a}$, where $\mathbf x(\tau)$, $0 \le \tau \le t$ is a path with $\mathbf{x}(0) = x$. At each time along the path, denote the direction of travel by $\theta(\tau)$. As discussed above, the travel velocity at the point $\mathbf x(\tau)$ and in the direction $\theta(\tau)$ is given by $V(s(\theta(\tau))\cdot\nabla E(\mathbf x(\tau))).$ Thus, the problem can be phrased as such: maximize the payoff function \begin{equation} P_{x,t}(\theta(\cdot)) =  \abs{\mathbf x(t) - a}\end{equation} among measurable functions $\theta: [0,t] \to [0,2\pi]$ and subject to the constraint \begin{equation} \begin{aligned} \dot{\mathbf{x}}(\tau) &=V(s(\theta(\tau))\cdot\nabla E(\mathbf x(\tau))) \theta(\tau), \,\,\,\, 0 \le \tau \le t \\ \mathbf x(0) &= x. \end{aligned} \end{equation} It is readily verified  that the Hamilton-Jacobi-Bellman equation associated with the value function $u(x,t) = \sup_{\theta} P_{x,t}(\theta(\cdot))$  for this control problem is precisely \ref{eq:levelSet} with the optimal path Hamiltonian \ref{eq:hamiltonian} and initial condition $\phi_0(x) = \abs{x-a}$. We then make the slight modification $\phi_0(x) = \abs{x-a}- \delta$ for small positive $\delta$ so that we may utilize the level set method to track optimal travel away from $a$ for every point on $\Gamma(0)$ simultaneously.

\subsection{Accounting for Uncertainty in the Starting Point} The above algorithm will compute a path for one who wishes to travel optimally throughout a region. We would like to incorporate some uncertainy into the model to account for a real world situation which law enforcement agents may encounter.Consider a scenario wherein a law enforcement agency has knowledge that environmental criminals (for example, poachers or illegal loggers) are operating within a protected region but can only identify the criminals' location with some uncertainty. Supposing that the criminals perpetrate a crime within the region and then travel to a known final destination, the law enforcement agency may want to predict which paths the criminals will take. 

\begin{figure}[b!]
\centering
\includegraphics[width=0.63\textwidth]{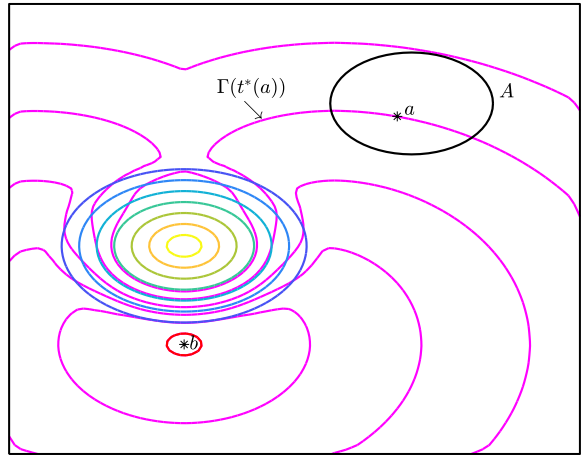}
\caption{If there is uncertainty in the location of the starting, we evolve level sets outward from $b$ until they cover $A$, recording optimal times for each $a \in A$ as we go.}
\label{withUncertainty}
\end{figure}

In this situation, assume that rather than a starting point $a$, we have a compact set $A$ of possible starting points along with a probability distribution from which one can sample elements of $A$. The algorithm described above requires a starting point which could be drawn from $A$ upon which one could calculate an optimal path to the end point $b$. However, we wish to calculate the optimal path to $b$ from each point in $A$ and according to our procedure, this will require solving \ref{eq:levelSet} for each point in $A$. Computationally, this would be very inefficient so instead, one can invert the problem: rather than starting from a point $a \in A$ and evolving level sets outward toward $b$, one should evolve level sets outward from $b$. As the level sets $\Gamma(t)$ evolve outward from $b$, they sweep through the region $A$ so that for each $a \in A$, we find a time $t^*(a)$ such that $a \in \Gamma(t^*(a))$. The computation can be stopped when $A$ is inside $\Gamma(t)$ and for each $a \in A$, we will have found the time $t^*(a)$ required to travel from $a$ to $b$.  Having done this, we can draw points $N$ points $\{a_i\}^N_{i=1}$ from $A$ and calculate the optimal paths using \ref{eq:ODE} with starting values $x(0) = a_i$, $p(0) = \nabla \phi(a_i, t^*(a_i))$. In this way, we can calculate optimal paths to $N$ points in $A$ with only one level set computation. 
 
As a minor note, walking velocity is maximal when one is walking on a slight decline. Thus reversing the direction of the level set evolution means we must also reverse our sense of slope since the walking direction is now opposite the direction of the level set evolution. Hence, we replace $E$ with $-E$. In doing so, when we evolve level sets outward from $b$, we are actually calculating optimal travel as if one was traveling inward toward $b$. Thus we can still compute the optimal path from $a$ to $b$ even though we use $b$ as the ``starting point'' for the level sets.

\section{Numerical Framework}\label{numerics} We discuss in detail the numerical methods that we use to simulate our model. The first obstacle is deciding how to calculate our Hamiltonian since this requires a maximization over $\theta \in [0,2\pi]$. If the velocity function $V$ is sufficiently simple, it may be possible to resolve this maximization explicitly using calculus. When this is not possible (as with our simulations), one can maximize $H$ discretely. That is, rather than maximize over $\theta \in [0,2\pi]$, we maximize $H_\theta(x,p)$ over the finite set $\theta \in \{ \tfrac{2\pi m}{M} \, : \, m = 1,\ldots, M\}.$ This causes some approximation error, but as long as $V(\nabla E(x) \cdot s(\theta))$ is continuous in $\theta$ for fixed $x$, this discrete maximization will tend to the exact supremum as $M \to \infty$.   

Next, one must decide how to solve \ref{eq:levelSet} numerically. There has been much research into efficient and accurate numerical methods for Hamilton-Jacobi equations \cite{AltonMitchell, kao2004lax, Osher1999, SzpiroDupuis}. Since these equations are (in general) nonlinear, naive differencing methods will not always work. Instead, we trade the Hamiltonian $H(x,\phi_x,\phi_y)$ for a numerical Hamiltonian $\hat H(x,\phi_x^+, \phi_x^-, \phi_y^+, \phi_y^-)$ which somehow averages the forward difference and backward difference approximations to $\phi_x$ and  $\phi_y$, represented here by $\phi_x^+,\phi_x^-$ and $\phi_y^+, \phi_y^-$ respectively.  We then advance the PDE via explicit time-stepping. Osher and Shu \cite{OsherShu} give several suggestions for different types of numerical Hamiltonians and describe methods for attaining higher order accuracy. For our purposes, we use the Godunov Hamiltonian given by \begin{equation}\hat H(x, \dphi{x}{+},\dphi{x}{-},\dphi{y}{+},\dphi{y}{-}) = \ext{u \in I(\dphi{x}{-},\dphi{x}{+})}\,\,\, \ext{v \in I(\dphi{y}{-},\dphi{y}{+})} H(x,u,v) \end{equation} where \begin{equation}I(a,b) = [\min(a,b), \max(a,b)]\end{equation} and \begin{equation} \ext{x \in I(a,b)} = \left \{\begin{matrix} \min_{a \le x \le b} & \text{if } a \le b, \\ \max_{b \le x \le a} & \text{if } a > b. \end{matrix} \right. \end{equation} These extrema are designed to take into account the direction in which information is flowing and, as a result, the Godunov Hamiltonian gives a fully upwind scheme. Again, we need to perform minimization or maximization computationally and again, in certain cases, these can be resolved explicitly (for example, if $H(x,u,v)$ is monotone in the arguments $(u,v)$), but this is not possible in our case, so we do this discretely. The Godunov Hamiltonian $\hat H$ gives a first-order approximation to the Hamiltonian $H$. Following Osher and Shu \cite{OsherShu}, we use second order essentially non-oscillatory approximations for the derivatives $\phi_x,\phi_y$ and second-order total variation diminishing Runge-Kutta time stepping to evolve the solution. In doing so, we have constructed a second order accurate scheme for \ref{eq:levelSet}. Finally, one can solve the optimal path ODE system \ref{eq:ODE} using any method one wishes. For relatively jagged elevation data $E$, the equation can become stiff, so it is recommended that one uses a stiff solver with accuracy which matches that of the numerical solution to \ref{eq:levelSet}. While this describes the basics of the numerical implementation, there are some minor adjustments required to obtain our results which we discuss in section \ref{implementation}.\\

\section{Implementation \& Results}\label{impRes} The model was implemented in MATLAB and in the succeeding section we discuss the results of the simulations and some issues which one may encounter. Before this, it remains to decide what elevation data to use and what form the velocity function takes. 

\subsection{Elevation \& Velocity}\label{env} 
\begin{figure}[htbp]
\centering
  \includegraphics[width=0.8\textwidth]{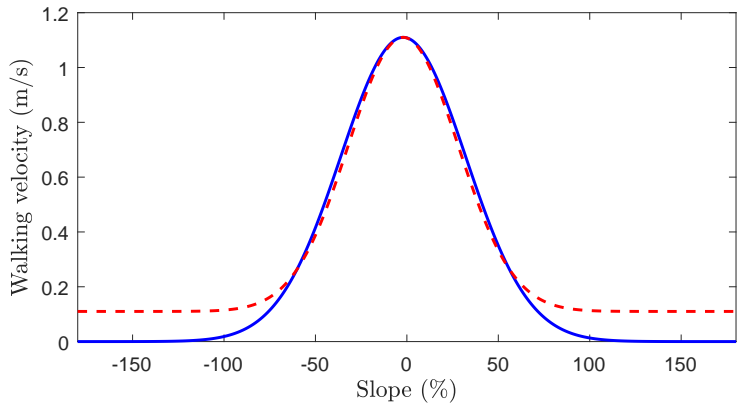} 
  \caption{Comparison of our velocity function (blue, solid) with Irmischer \& Clarke's velocity function (red, dotted).} 
  \label{velfunc}
\end{figure}

For the velocity function, we use a slight modification of the function Irmischer and Clarke \cite{IC2017}. Irmischer and Clarke analyze human walking speed data and suggest the function \begin{equation}V_{IC}(S) = 0.11 + \exp\left(-\frac{(100S + 2)^2}{1800}\right)\end{equation} where $S = \frac{\text{rise}}{\text{run}}$. However, in their paper, they only considered slopes up to $45^\circ$ (grades up to 100\%), and their function is bounded below by $0.11$.  This is a good starting point, but we would like to consider slopes much higher than $45^\circ$ where walking speed may become very small. Accordingly, using a slightly different ansatz and fitting the denominator in the exponential, we have arrived at our own velocity function which approximates the Irmischer and Clarke function for small slopes but which decays to zero for more extreme slopes: \begin{equation}V(S) = 1.11 \exp\left(-\frac{(100S + 2)^2}{2345}\right). \end{equation}While this function is never exactly zero, it is no longer bounded from below by any positive number. It bears mentioning that the exact form of the velocity function is not terribly important for the model so long as the function $V(S)$ that we choose has $\max V(S) \approx 1$, is fairly near the maximum for all $S \approx 0$ and is nearly zero for $\abs S$ large. Our velocity function is plotted against the Irmischer-Clarke velocity function in Figure~\ref{velfunc}.

To test our code, we first ran simulations with synthetic (and very simple) elevation data. This allows us to guage whether our model aligns with our intuition. When we were confident that our model and numerical methods were correct, we were able to download real elevation data from the United States Geological Survey and run simulations in a real national park. For our simulations, we chose Yosemite National Park and we ran optimal path simulations in the direct vicinity of the mountain El Capitan. Specifically, we use data spanning longitude $119^\circ$W - $120^\circ$W and lattitude $38^\circ$N - $39^\circ$N with 1/3 arcsecond resolution which we obtained from the USGS National Map Viewer. The data was processed and re-formatted using QGIS \cite{QGIS} and imported to MATLAB using TopoToolbox \cite{topotoolbox2,topotoolbox1}. \\

\subsection{Results} \label{results}
\begin{figure}[htbp]
\centering
\includegraphics[width=0.63\textwidth]{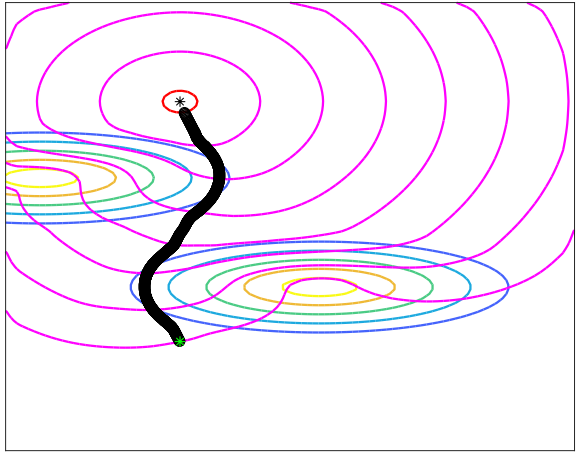}
\caption{Optimal path winding around two mountains (toy problem).}
\label{twomountains}
\end{figure}

As in Figure~\ref{key1} about, in  the following images, the starting point $a$ is represented by the black asterisk surrounded with a red circle which denotes the starting contour $\Gamma(0)$. Next, the magenta contours represent several steps in the evolution of the contours $\Gamma(t)$. The green asterisk represents the point $b$ and the thick black line represents the optimal path from $a$ to $b$. The elevation contours are plotted in colors ranging from blue representing low elevation to yellow representing high elevation. In our first simulations, we place mountains in certain areas and our intuition tells us that the optimal path should likely bend around the mountains since it would require too much effort to climb up the mountain. Our code does indeed predict this; see Figure~\ref{twomountains}.

\begin{figure}[b!]
\centering
\includegraphics[width=0.63\textwidth]{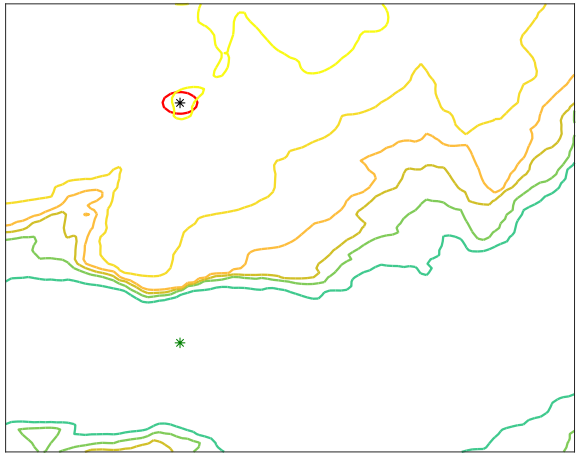}
\caption{Elevation profile of El Capitan.}
\label{elcapelev}
\end{figure}

\begin{figure}[htbp]
\centering
\begin{subfigure}[t]{0.45\textwidth}
\centering
\includegraphics[width=0.95\textwidth]{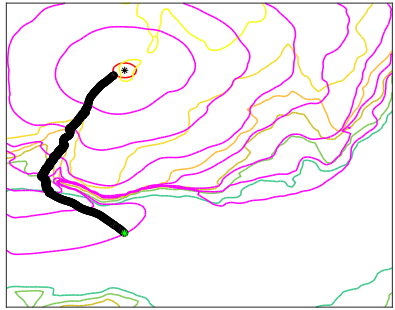}
\caption{Optimal path down the western slope.}
\label{elcap1}
\end{subfigure}
\qquad
\begin{subfigure}[t]{0.45\textwidth}
\centering
\includegraphics[width=0.95\textwidth]{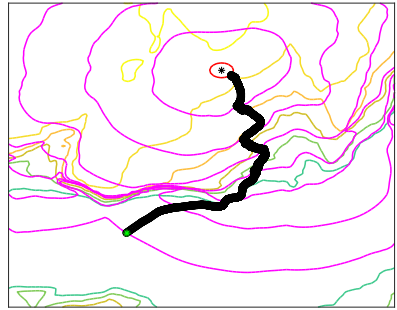}
\caption{Optimal path down the eastern slope.}
\label{elcap2}
\end{subfigure}
\label{elcap}
\caption{Optimal paths down from El Capitan avoid the steep cliff.}
\end{figure}

Next, we use actual elevation data from the area surrounding El Capitan. Before showing the result of the simulation, we show the elevation profile and the starting and ending points in Figure~\ref{elcapelev}. Note that directly above the endpoint, there is a very steep cliff face which should be nearly impossible to traverse. Thus we would expect the optimal path to travel to the east or west, descending down a gully rather than a cliff.  Indeed, this is shown to happen in Figures~\ref{elcap1},\ref{elcap2}, wherein the path travels down the eastern or western slope depending on the location of the initial point.

\begin{figure}[htbp]
\centering
\begin{subfigure}[t]{0.45\textwidth}
\centering
\includegraphics[width=0.95\textwidth]{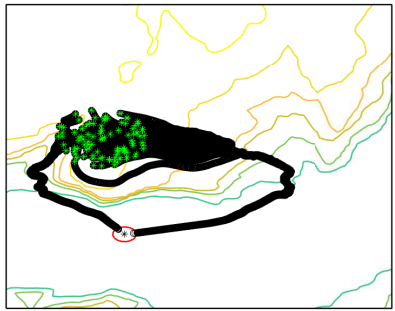}
\caption{100 optimal paths traveling down from the summit of El Capitan.}
\label{elCapUncertainty}
\end{subfigure}
\qquad
\begin{subfigure}[t]{0.45\textwidth}
\centering
\includegraphics[width=0.95\textwidth]{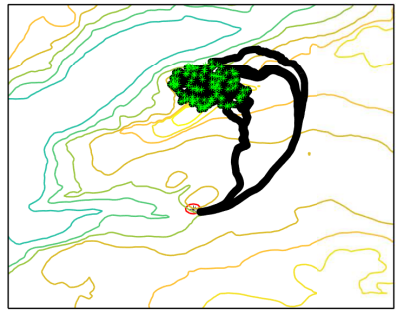}
\caption{100 optimal paths traveling down from the summit of Half Dome.}
\label{halfDomeUncertainty}
\end{subfigure}
\label{Uncertainty2}
\caption{Calculation optimal paths accounting for uncertainty in the initial location.}
\end{figure}

Next, we ran our algorithm which accounts for uncertainty in the location of the starting point. Before we display our results, we remind the reader of the distinction here. In all of these above results, we are calculating one optimal path from the point $a$ to the point $b$. Now we wish to calculate several optimal paths from the region $A$ to the point $b$. Whereas previously we evolved level sets outward from point $a$ until they reach the point $b$, now we evolve level sets outward from $b$ until they subsume the region $A$ and record the optimal travel time for each $a \in A$ as the level sets sweep through the region. This is shown in Figures~\ref{withUncertainty}. We ran our algorithm in two different areas within Yosemite. We let $A$ be a circle near the peak of El Capitan and calculated the optimal path down the mountain from 100 random points drawn uniformly from $A$. We then did the same thing but using the elevation profile of Half Dome, another peak in Yosemite National Park. The results are pictured in Figures~\ref{elCapUncertainty},\ref{halfDomeUncertainty}. 

% \begin{figure}[htbp]
%\centering
%\begin{subfigure}[t]{0.5\textwidth}
%\centering
%\includegraphics[width=0.95\textwidth]{withoutUncertainty.eps}
%\caption{If there is no uncertainty, we evolve level sets outward from $a$ until they reach $b$.}
%\label{withoutUncertainty}
%\end{subfigure}~
%\begin{subfigure}[t]{0.5\textwidth}
%\centering
%\includegraphics[width=0.95\textwidth]{withUncertainty.eps}
%\caption{If there is uncertainty, we evolve level sets outward from $b$ until they cover $A$, recording optimal times for each $a \in A$ as we go.}
%\label{withUncertainty}
%\end{subfigure}
%\label{Uncertainty}
%\caption{Difference in the algorithm if there is no uncertainty, (a), or is uncertainty, (b), in the location of the starting point.}
%\end{figure}  

  \begin{figure}[htbp]
\centering
\begin{subfigure}[t]{0.45\textwidth}
\centering
\includegraphics[width=0.95\textwidth]{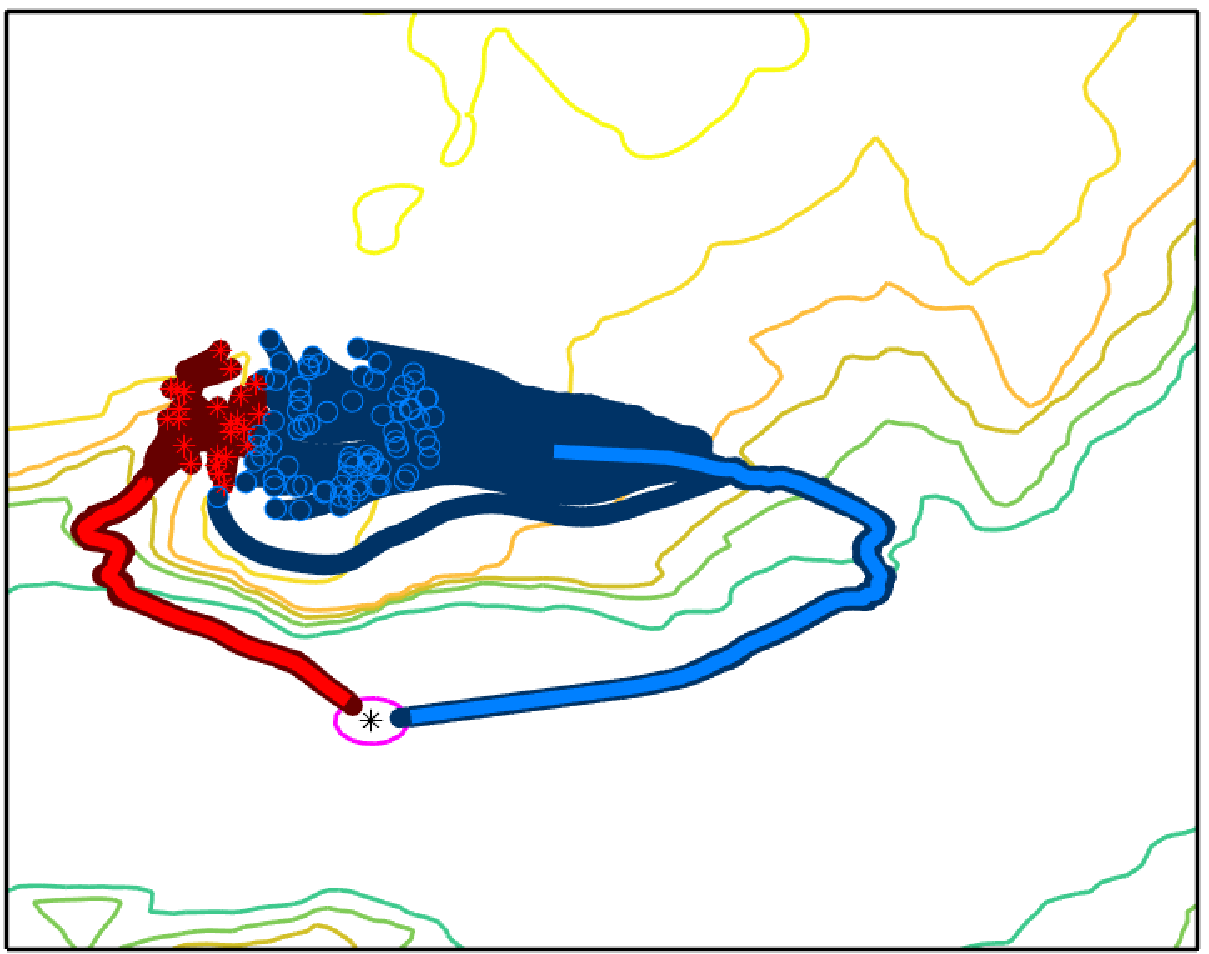}
\caption{Clustering the paths down from El Capitan into two collections.}
\label{elCapPathClusters}
\end{subfigure}
\qquad
\begin{subfigure}[t]{0.45\textwidth}
\centering
\includegraphics[width=0.95\textwidth]{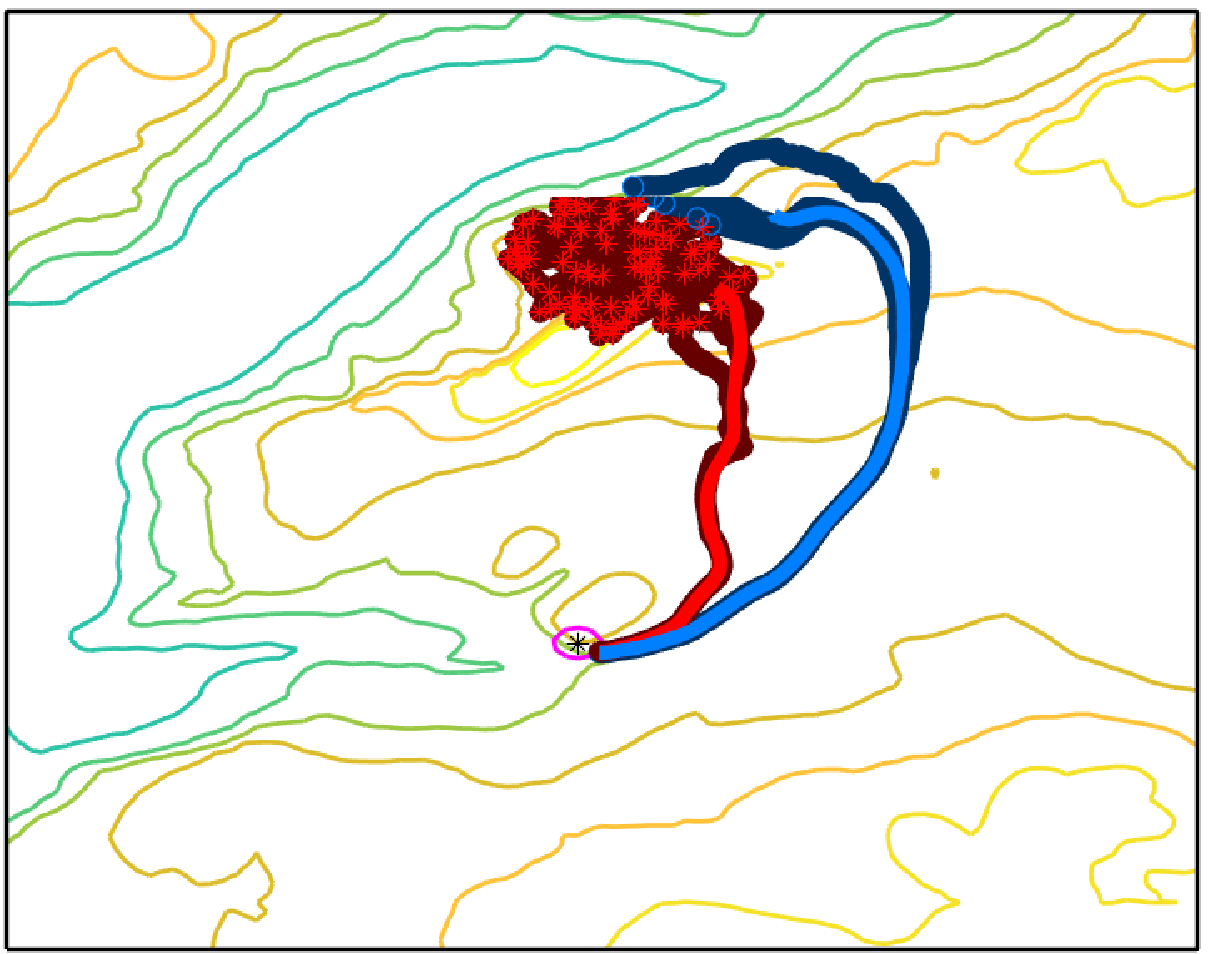}
\caption{Clustering the paths down from Half Dome into two collections.}
\label{halfDomePathClusters}
\end{subfigure}
\label{PathClusters}
\caption{Clustering can help us identify which paths are morally the same. The bright blue path is the representative path for the blue points and the bright red path is the representative path for the red points.}
\end{figure}

Note that in both cases, while there are 100 randomly chosen starting positions, all of the paths eventually conform to one of very few routes. We seek to quantify this similarity between some paths. Having calculated several paths $\{ P_i\}_{i=1}^N$ starting from different locations, one can re-parametrize so that for each $i$, $P_i: [0,1] \to \R^2$ and $P(0) = b$. Then it is easy to define a metric to judge whether two paths lie nearby each other: for two paths $P,Q$, define $d(P,Q) = \int^{0.8}_0 \norm{P(t) - Q(t)} \mathrm{d}t$. The reason for this $0.8$ in the upper limit is so that we do not penalize paths for starting from different points. With this metric, one can evaluate the pairwise distances between our paths, $\{d(P_i,P_j)\}_{i,j=1}^N$. Now, using basic clustering algorithms, one can categorize the paths into collections which are morally the same, in the sense that they eventually collapse onto the same route. We performed this clustering for the above two examples. Specifically, we used $k$-means clustering with $k = 2$ clusters in each case, though other clustering methods could be used. The results are included in Figures~\ref{elCapPathClusters},\ref{halfDomePathClusters}, where the first cluster of paths is depicted in red and the second in blue. Here, rwe have plotted each of the 100 paths as well as the mean path for each cluster. Thus any initial point which is marked with a blue circle has a corresponding optimal path which eventually closely resembles the bright blue path and any initial point which is marked with a red asterisk has a corresponding optimal path which eventually closely resembles the bright red path. Returning to our original motivation, these graphics could be of great interest to law enforcement agencies who are tracking criminal movement. For example, in the case of El Capitan (Figure~\ref{elCapPathClusters}), we observe that $24\%$ of the paths travel down the eastern slope while $76\%$ travel down the eastern slope. This may suggest to law enforcement that they should patrol the eastern slope with roughly three times the resources which they devote to the western slope.

\subsection{Implementation Notes}\label{implementation}

\begin{figure}[htbp]
\centering
\begin{subfigure}[t]{0.45\textwidth}
\centering
\includegraphics[width=0.95\textwidth]{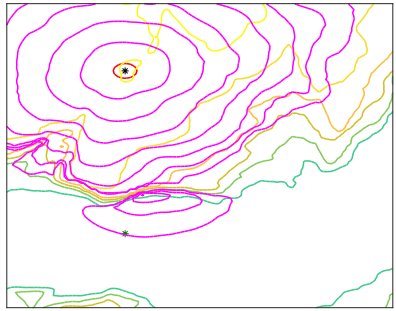}
\caption{Level sets without re-distancing ``jump over" the cliff.}
\label{noredist}
\end{subfigure}
\qquad
\begin{subfigure}[t]{0.45\textwidth}
\centering
\includegraphics[width=0.95\textwidth]{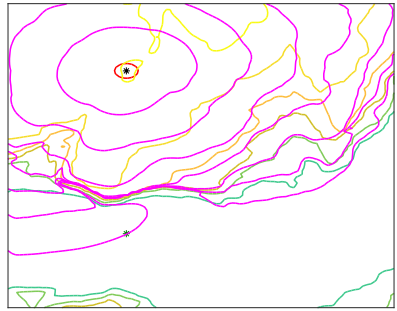}
\caption{Level sets with re-distancing wrap around the cliff.}
\label{withredist}
\end{subfigure}
\caption{Level sets without, (a), and with, (b), re-distancing.}
\label{redist}
\end{figure}

\begin{figure}[htbp]
\centering
\begin{subfigure}[t]{0.45\textwidth}
\centering
\includegraphics[width=.95\textwidth]{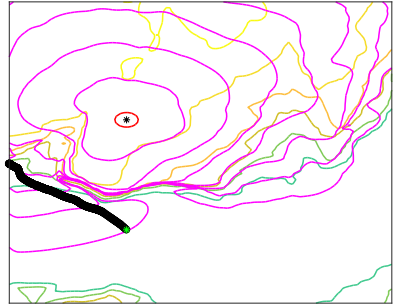}
\caption{Optimal path without re-initialization veers of the map.}
\label{noreinit}
\end{subfigure}
\qquad
\begin{subfigure}[t]{0.45\textwidth}
\centering
\includegraphics[width=.95\textwidth]{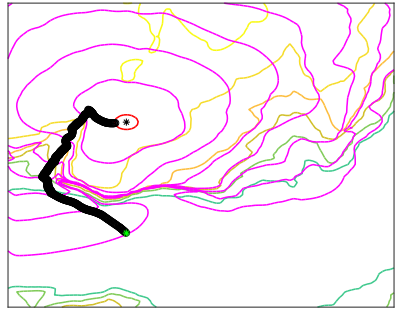}
\caption{Optimal path with re-initialization finds its target.}
\label{withreinit}
\end{subfigure} 
\caption{Optimal paths without, (a), and with, (b), re-initilialization.}
\label{reinit}
\end{figure}

\begin{figure}[htbp]
\centering
\begin{subfigure}[t]{0.45\textwidth}
\centering
\includegraphics[width=.95\textwidth]{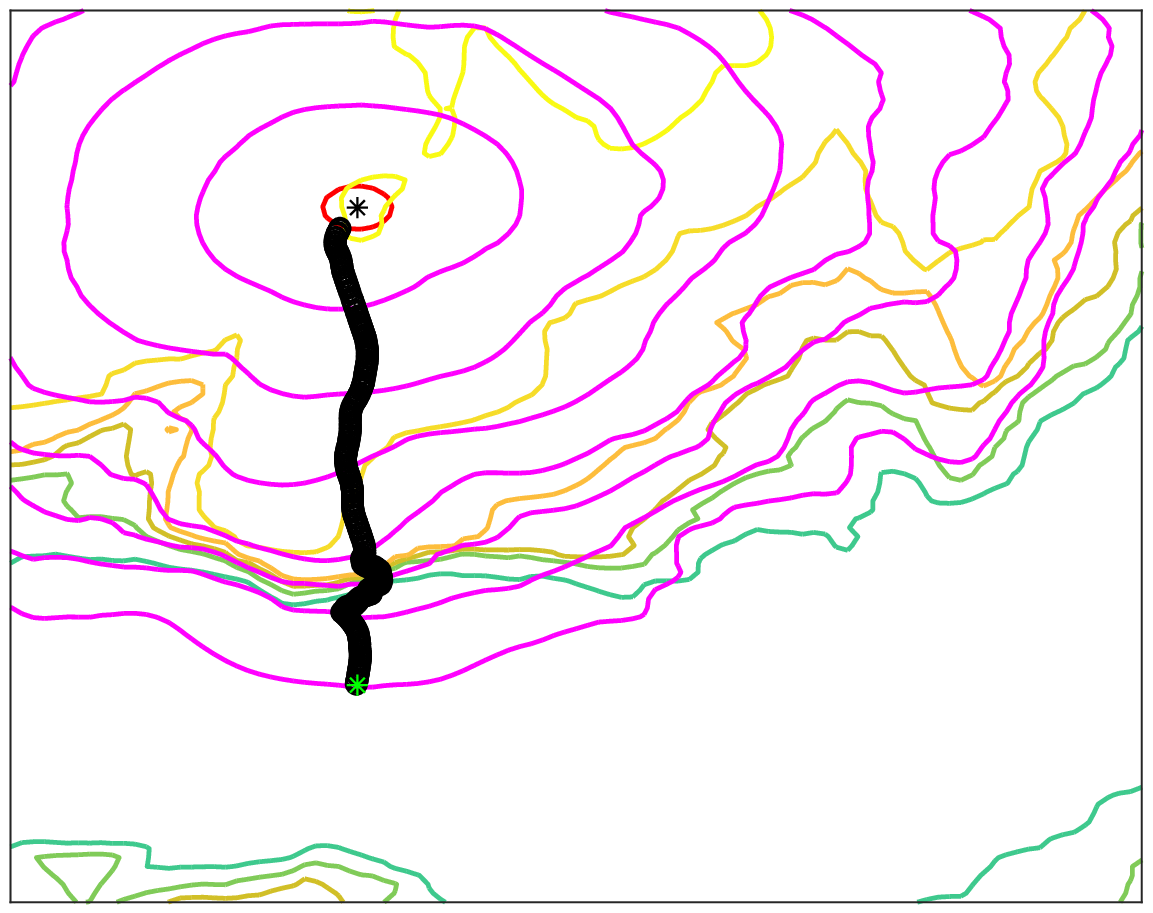}
\caption{Optimal path without high-slope penalization zig-zags up the cliff.}
\label{novelcorrect}
\end{subfigure}
\qquad
\begin{subfigure}[t]{0.45\textwidth}
\centering
\includegraphics[width=.95\textwidth]{ElCapOpPath1.eps}
\caption{Optimal path with high-slope penalization avoids the cliff.}
\label{withvelcorrect}
\end{subfigure} 
\caption{Optimal paths without, (a), and with, (b), the high-slope penalization.}
\label{slopepen}
\end{figure}

There are a few specific issues that arise when implementing the model numerically. We discuss three such issues (and their resolutions) and demonstrate their effects in Figures~\ref{noredist}-\ref{withvelcorrect}. First, note that the initial function $\phi(x,0)$ gives precisely the signed distance from $x$ to $\Gamma(0)$; that is $$\phi(x,0) =\text{dist}(x,\Gamma(0)) \defeq  \left\{\begin{matrix} \inf_{y \in \Gamma(0)} \abs{x-y}, & x \text{ inside } \Gamma_0, \\ -\inf_{y \in\Gamma(0)} \abs{x-y}, & x \text{ outside } \Gamma(0). \end{matrix} \right. $$ As the level sets evolve, there is some distortion so that for $t > 0$ we no longer have $\phi(x,t) = \text{dist}(x,\Gamma(t))$. This distortion happens when $\abs{\nabla \phi}$ becomes too large or too small near the zero level contour $\Gamma(t)$ and can cause the level set results to become unreliable. We can fix this by occasionally replacing $\phi$ with the signed distance function to $\Gamma(t)$. That is, we occasionally halt the time integration, find the current zero level contour $\Gamma(t)$, reset $\phi(x,t) = \text{dist}(x,\Gamma_t)$ and continue. This process is known as \emph{re-distancing} and is discussed in several papers \cite{Osher2012, redistancing, SussmanFatemi}. Figure~\ref{redist} shows the effect of re-distancing.

Next, as mentioned before, the system \ref{eq:ODE} used to find the optimal path becomes very stiff when non-smooth elevation profiles are used. Even when using a stiff solver, the results were unreliable in that the value of $p(t)$ corresponding to a location $x(t)$ was straying far from the theoretically correct value $\nabla \phi_x(x(t),t)$. This was causing the ``optimal path" that our code found to be wildly inaccurate, often times not even connecting $b$ to $a$, opting instead to wander off in some seemingly random direction. To fix this, we do something similar to the above fix: we occasionally stop the time integration, re-initialize $p(t) = \nabla \phi(x(t),t)$ and then restart the time stepping. We refer to this as \emph{re-initialization} and the effect is shown in Figure~\ref{reinit}. 

Finally, there is still one shortcoming of our Hamiltonian with respect to directional movement: the slope in the direction of travel and its orthogonal are completely decoupled. For example, consider a situation where there is a steep cliff face in the north-south direction while the slope in the east-west direction is very mild. Our model would allow an individual to walk east-west in this situation even though they may be standing on an prohibitively steep slope. To fix this problem, we add a penalty for walking in locations where the maximum slope in \emph{any} direction is very large. This is as simple as multiplying our Hamiltonian by a pre-factor which is approximately $1$ for low slopes and approximately zero for high slopes. We have chosen the penalization function $P(S) = \tfrac 1 2 - \tfrac 1 2\tanh(S - 1)$ where $S = \frac{\text{rise}}{\text{run}}$ is the slope. Thus we actually solve the Hamilton-Jacobi-Bellman equation with Hamiltonian $P(\abs{\nabla E(x)})H(x,p)$ where $H(x,p)$ is the optimal path Hamiltonian.

\section{Conclusion} \label{conc} We have presented a method for resolving optimal walking paths given terrain data. The key element of the model is a generalization of the level set equation. By representing the direction of travel for the level sets with a control variable, we constructed a Hamiltonian whose corresponding level set equation models optimal travel. Using this, we described a simple algorithm for calculating the optimal walking path between a starting and ending point which consists of numerically solving a Hamilton-Jacobi-Bellman (HJB) equation and then a system of ordinary differential equations. Further, we suggest a method for incorporating uncertainty into the location of the starting point: by modifying the algorithm slightly, we can compute several optimal paths while only solving one HJB equation. We then suggested numerical methods for simulating the HJB equation. We used Godunov's numerical Hamiltonian with second order essentially non-oscillatory finite difference approximations for spatial derivatives and second order total variation diminishing time integration. We also suggested modifications to the numerical methods which avoid common pitfalls which one may encounter. To test our algorithm, we simulated our model first using artificial elevation data and then using the actual elevation data in certain regions of Yosemite National Park. In both cases, results aligned very well with our physical intuition. Finally, we sampled several different starting locations and calculated optimal paths which travel down from the summits of El Capitan and Half Dome and noticed that in both cases, the paths can be naturally clustered into collections of paths which follow the same basic route. We performed $k$-means clustering to separate the paths into such collections. Such clustering could suggest simple yet effective patrol strategies for a law enforcement agency tasked with patrolling nationally protected areas. 

\section{Acknowledgements} Author Yat Tin Chow is supported by the National Science Foundation (DMS-1737770 and ECCS-1462398). The remaining authors are funded by the National Geospatial-Intelligence Agency's Academic Research Program (Award No. HM0210-14-1-0003, Project Title: ``Sparsity Models for Spatiotemporal Analysis and Modeling of Human Activity and Social Networks in a Geographic Context." Approved for public release,
18-648.

\bibliographystyle{plainnat}
\bibliography{bibliography}
\end{document}